\documentclass{amsart}%
\usepackage{amsfonts}%
\usepackage{amsmath}%
\usepackage{amssymb}%
\usepackage{graphicx}
\newtheorem{theorem}{Theorem}
\theoremstyle{plain}

\newtheorem{corollary}{Corollary}

\newtheorem{definition}{Definition}

\newtheorem{problem}{Problem}
\newtheorem{proposition}{Proposition}
\newtheorem{remark}{Remark}

\numberwithin{equation}{section}
\begin{document}
\title[Compactifications of Riemannian manifolds]{Compactifications of Complete Riemannian manifolds and Their Applications}
\author{Xiaodong Wang}
\address{Department of Mathematics\\
Michigan State University\\
East Lansing, MI\ 48824}
\email{xwang@math.msu.edu}
\thanks{The author was partially supported by NSF grant DMS-0905904.}
\dedicatory{Dedicated to Professor Richard Schoen in honor of his 60th birthday}\maketitle

\section{\bigskip Introduction}

To study a noncompact Riemannian manifold, it is often useful to find a
compactification or attach a boundary. For example, in hyperbolic geometry a
lot of investigation is carried out on the sphere at infinity. An eminent
illustration is Mostow's proof of his rigidity theorem for hyperbolic
manifolds \cite{Mo}. More generally, if $\widetilde{M}$ is simply connected
and nonpositively curved, one can compactify it by equivalent geodesic rays
and the boundary is a topological sphere, called the geometric boundary. This
compactification was first introduced in \cite{EO} and has been indispensable
in the study of negatively curved manifolds. If $\widetilde{M}$ is not
nonpositively curved, then the geometric compactification does not work in
general. But there are other compactifications which are useful for various
studies. In this short survey, we will discuss some of these compactifcations
and the relationships among them. Our discussion will focus on general
Riemannian manifolds and therefore we ignore the large literature on
compactifications of symmetric spaces (see the book \cite{GJL}).

We first discuss the geometric compactification for Cartan-Hadamard manifolds
and Gromov hyperbolic spaces in Section 2. In Section 3 we discuss the Martin
compactification. In Section 4 we discuss the Busemann compactification. In
the last section, we discuss how these compactifications are used. In
particular, we consider certain invariants defined on the Martin boundary and
prove a comparison inequality using a method of Besson, Courtois and Gallot.
It should be noted that when the author showed this inequality to Fran\c{c}ois
Ledrappier he was informed that it had been known to Besson, Courtois and
Gallot (unpublished).

\bigskip

\textit{Acknowledgement:} It is my great pleasure to dedicate this article to
Professor Rick Schoen on the occasion of his 60th birthday. Professor Schoen's
many fundamental contributions to geometry are well-known and greatly admired.
I had the good fortune to be one of his students and his influence on my
mathematical career, through his teaching, his work and his example, has been
tremendous. I wish him good health and more great theorems in the years to come.

In preparing this article, I have benefited from several stimulating
discussions with Fran\c{c}ois Ledrappier. I want to thank him warmly for
teaching me a lot of things. I would also like to thank Jianguo Cao for
helpful discussion on the work of Ancona \cite{A} and for bringing to my
attention his paper \cite{Cao}.

\section{The geometric compactification}

The most familiar compactification is the geometric compactification first
introduced by Eberlein and O'Neill \cite{EO} for Cartan-Hadamard manifolds.
Let $\widetilde{M}^{n}$ be a Cartan-Hadamard manifold. We can compactify
$\widetilde{M}$ using geodesic rays. More precisely, two geodesic rays
$\gamma_{1}$ and $\gamma_{2}$ are said to be equivalent, if $d\left(
\gamma_{1}\left(  t\right)  ,\gamma_{2}\left(  t\right)  \right)  $ is bounded
for $t\in\lbrack0,\infty)$. The set of equivalence classes, denoted by
$\widetilde{M}\left(  \infty\right)  $, is called the geometric boundary and
can be naturally identified with the unit sphere $\mathbb{S}^{n-1}$ if we fix
a base point. We then obtain a compactification $\widetilde{M}^{\ast
}=\widetilde{M}\sqcup\widetilde{M}\left(  \infty\right)  $ that is
homeomorphic to the closed unit $n$-ball with the natural "cone" topology. If
the sectional curvature satisfies $-b^{2}\leq K_{\widetilde{M}}\leq-a^{2} $,
where $a,b>0$, it is proved by Anderson and Schoen that $\widetilde{M}^{\ast}$
has a $C^{\alpha}$-structure, where $\alpha=a/b$. For details, cf. \cite{E,
SY}.

The same compactification works for the so called Gromov hyperbolic spaces. We
first recall one of several equivalent definitions of Gromov hyperbolic spaces.

\begin{definition}
A complete geodesic metric space $\left(  X,d\right)  $ is called Gromov
hyperbolic if for some $\delta>0$ s.t. for all points $o,x,y,z\in X$%
\[
\left(  x\cdot y\right)  _{o}\geq\min\left\{  \left(  x\cdot z\right)
_{o},\left(  y\cdot z\right)  _{o}\right\}  -\delta\text{.}%
\]
where we use Gromov products, e.g. $\left(  x\cdot y\right)  _{o}=\frac{1}%
{2}\left(  d\left(  o,x\right)  +d\left(  o,y\right)  -d\left(  x,y\right)
\right)  $.
\end{definition}

We make the following remarks

\begin{remark}
In the definition one can take $o$ to be fixed.
\end{remark}

\begin{remark}
A Cartan-Hadamard manifold $\widetilde{M}$ is Gromov hyperbolic if the
sectional curvature has a negative upper bound.
\end{remark}

This concept was introduced by Gromov \cite{G}. For detailed study of Gromov
hyperbolic spaces, see the excellent books \cite{BH, GH, O}. It suffices to
say that the definition captures the global features of the geometry of a
complete simply connected manifold of negative curvature. It is very robust as
illustrated by the following remarkable fact.

\begin{theorem}
Let $X$ and $Y$ be geodesic spaces. Suppose $f:X\rightarrow Y$ is a
quasi-isometry, i.e. there are $L,\varepsilon>0$ s.t. for any $x_{1},x_{2}\in
X$
\[
L^{-1}d\left(  x_{1},x_{2}\right)  -\varepsilon\leq d\left(  f\left(
x_{1}\right)  ,f\left(  x_{1}\right)  \right)  \leq Ld\left(  x_{1}%
,x_{2}\right)  .
\]
If $Y$ is Gromov hyperbolic, then so is $X$.
\end{theorem}

We will further assume that $X$ is proper. Then we can define the geometric
boundary $X\left(  \infty\right)  $ for a Gromov hyperbolic spaces in the same
way s.t. $\overline{X}=X\sqcup X\left(  \infty\right)  $ with a natural
topology is a compact metrizable space. Moreover $X\left(  \infty\right)  $
has a canonical quasi-conformal structure.

\begin{theorem}
Let $X$ and $Y$ be proper Gromov hyperbolic spaces. If $f:X\rightarrow Y$ is a
quasi-isometry, then $f$ extends to a homeomorphism $\overline{f}:X\left(
\infty\right)  \rightarrow Y\left(  \infty\right)  $.
\end{theorem}

In fact, the boundary map is furthermore a quasi-conformal map. The boundary
map can be described as follows: given the equivalence class $\xi\in\partial
X$ of a geodesic ray $\gamma:[0,\infty)\rightarrow X$, $f\circ\gamma$ is a
quasi-geodesic in $Y$ and hence has a well defined end point $f\circ
\gamma\left(  \infty\right)  \in Y\left(  \infty\right)  $ which is defined to
be $\overline{f}\left(  \xi\right)  $.

\section{The Martin compactification}

We can also compactify $\widetilde{M}$ using all positive harmonic functions.
The vector space $\mathcal{H}\left(  \widetilde{M}\right)  $ of harmonic
functions with seminorms
\[
\left\Vert u\right\Vert _{K}=\sup_{K}\left\vert u\left(  x\right)  \right\vert
,K\subset\widetilde{M}\text{ compact}%
\]
is a Frechet space. Let $\mathcal{K}_{o}=\{u\in\mathcal{H}_{+}\left(
\widetilde{M}\right)  :u\left(  o\right)  =1\}$. It is a convex and compact
set in $\mathcal{H}\left(  \widetilde{M}\right)  $. We assume that
$\widetilde{M}$ is nonparabolic and $G(x,y)$ is the minimal positive Green's
function. Define the Martin kernel%
\[
k\left(  x,y\right)  =\frac{G(x,y)}{G(o,y)}.
\]
A sequence \thinspace$y_{i}\rightarrow\infty$ is called a Martin sequence if
$\lim_{i\rightarrow\infty}k\left(  x,y_{i}\right)  $ converges to a harmonic
function. By Harnack inequality and the elliptic theory every sequence
$y_{i}\rightarrow\infty$ has a Martin subsequence. Two Martin sequences are
called equivalent if they have the same harmonic function as limit. The
collection of all such equivalence classes is called the Martin boundary and
will be denoted by $\partial_{\Delta}\widetilde{M}$. It is easy to see that
$\partial_{\Delta}\widetilde{M}\subset\mathcal{K}_{o}$ is a compact set. The
Martin compactification is defined to be%
\[
\widehat{M}=\widetilde{M}\sqcup\partial\widetilde{M}%
\]
with a natural topology that makes it a compact metrizable space. An excellent
reference on Martin compactification is Ancona \cite{A}.

\begin{definition}
A harmonic function $h>0$ on $\widetilde{M}$ is called minimal if any
nonnegative harmonic function $\leq h$ is proportional to $h$.
\end{definition}

\begin{remark}
If $h\left(  o\right)  =1$, then $h$ is minimal iff $h$ is an extremal point
of $\mathcal{K}_{o}$.
\end{remark}

It is proved that all minimal harmonic function $h$ with $h\left(  o\right)
=1 $ belong to $\partial_{\Delta}\widetilde{M}\,$. Therefore we can introduce
the following

\begin{definition}
The minimal Martin boundary of $M$ is
\[
\partial^{\ast}\widetilde{M}=\{h\in\mathcal{K}_{o}:h\text{ is minimal}\}.
\]

\end{definition}

Moreover $\partial^{\ast}\widetilde{M}\subset\partial\widetilde{M}$ is at
least a Borel subset (cf. \cite{A}). According to a theorem of Choquet
(\cite{A}), for any positive harmonic function $h$ there is a unique Borel
measure $\mu^{h}$ on $\partial^{\ast}\widetilde{M}$ such that%
\[
h\left(  x\right)  =\int_{\partial^{\ast}\widetilde{M}}\xi\left(  x\right)
d\mu^{h}\left(  \xi\right)  .
\]
Let $\nu$ be the measure corresponding to the harmonic function $1$. Thus%
\begin{equation}
1=\int_{\partial^{\ast}\widetilde{M}}\xi\left(  x\right)  d\nu\left(
\xi\right)  .\label{defnu}%
\end{equation}
The family of probability measures $\left\{  \nu^{x}:x\in\widetilde
{M}\right\}  $ with $\nu^{x}=$ $\xi\left(  x\right)  \nu$ are called the
harmonic measures. For $f\in L^{\infty}\left(  \partial^{\ast}\widetilde
{M}\right)  $ we get a bounded harmonic function%
\[
H_{f}\left(  x\right)  =\int_{\partial^{\ast}\widetilde{M}}f\left(
\xi\right)  \xi\left(  x\right)  d\nu\left(  \xi\right)  .
\]
This defines an isomorphism between $L^{\infty}\left(  \partial^{\ast
}\widetilde{M}\right)  $ and the space of bounded harmonic functions on
$\widetilde{M}$.

The study of the Martin compactification is closely related to the study of
Brownian motion on $\widetilde{M}$. For simplicity we further add a mild
condition that the Ricci curvature is bounded from below to ensure stochastic
completeness. Therefore we have the sample space $\Omega\left(  \widetilde
{M}\right)  =C\left(  \mathbb{R}^{+},\widetilde{M}\right)  $ with a family of
probability measure $\left\{  \mathbb{P}^{x}:x\in\widetilde{M}\right\}  $ s.t.
for any $0<t_{1}<\cdots<t_{k}$ and open sets $U_{1},\cdots,U_{k}$
\begin{align*}
& \mathbb{P}^{x}\left\{  \omega\in\Omega\left(  M\right)  :\omega\left(
t_{1}\right)  \in U_{1},\cdots\omega\left(  t_{k}\right)  \in U_{k}\right\} \\
& =\int_{U_{1}\times\cdots\times U_{k}}p_{t_{1}}\left(  x,y_{1}\right)
p_{t_{2}-t_{1}}\left(  y_{1},y_{2}\right)  \cdots p_{t_{k}-t_{k-1}}\left(
y_{k-1},y_{k}\right)  dy_{1}\times\cdots dy_{k}.
\end{align*}
Here $p_{t}\left(  x,y\right)  $ is the heat kernel on $\widetilde{M}$. For
each $t\geq0$ we have a random variable $X_{t}:$ $\Omega\left(  \widetilde
{M}\right)  \rightarrow\widetilde{M}$ which is simply the position at $t$,
i.e. $X_{t}\left(  \omega\right)  =\omega\left(  t\right)  $. It is an
intriguing and important problem to understand the asymptotic behavior for
$\omega\left(  t\right)  $ as $t\rightarrow\infty$. The answer is closely
related to the Martin boundary.

\begin{theorem}
\begin{enumerate}
\item For any $x\in\widetilde{M}$ and for $\mathbb{P}^{x}$-a.e. $\omega
\in\Omega\left(  M\right)  $, $X_{t}\left(  \omega\right)  $ admits a limit
$X_{\infty}\left(  \omega\right)  \in\partial^{\ast}\widetilde{M}$ as
$t\rightarrow\infty$, i.e.%
\[
\lim_{t\rightarrow\infty}k\left(  x,\omega\left(  t\right)  \right)
\]
exists and is a minimal harmonic function.

\item Under $\mathbb{P}^{x}$, the distribution of $X_{\infty}$ is $\nu^{x}$,
i.e. $\left(  X_{\infty}\right)  _{\ast}\mathbb{P}^{x}=\nu^{x}$.
\end{enumerate}
\end{theorem}

For detailed discussion, see Ancona \cite{A}.

For a Cartan-Hadamard manifold $\widetilde{M}$ with sectional curvature
bounded between two negative constants, Anderson and Schoen \cite{AS} proved
that the Martin boundary is homeomorphic to the geometric boundary.

\begin{theorem}
\bigskip Suppose that $\widetilde{M}$ is a Cartan-Hadamard manifold with whose
sectional curvature satisfies $-b^{2}\leq K\leq-a^{2}<0$. Then there exists a
natural homeomorphism $\Phi:\partial\widetilde{M}\rightarrow\widetilde
{M}\left(  \infty\right)  $ between the Martin boundary and the geometric
boundary. Moreover, $\Phi^{-1}$ is H\"{o}lder continuous.
\end{theorem}

From the proof, it is also clear that $\partial\widetilde{M}=\partial^{\ast
}\widetilde{M}$ in this case.

This theorem was generalized by Ancona who proved

\begin{theorem}
\label{Ancona}(Ancona \cite[Theorem 6.2]{A}) Suppose that $\widetilde{M}$ is
Gromov hyperbolic and $\lambda_{0}\left(  \widetilde{M}\right)  >0$. Then the
Martin boundary is homeomorphic to the geometric boundary. Moreover
$\partial\widetilde{M}=\partial^{\ast}\widetilde{M}$.
\end{theorem}

In the statement, $\lambda_{0}\left(  \widetilde{M}\right)  $ is the bottom of
the $L^{2}$ spectrum of $\widetilde{M}$, i.e.
\[
\lambda_{0}\left(  \widetilde{M}\right)  =\inf\frac{\int_{\widetilde{M}%
}\left\vert \nabla u\right\vert ^{2}}{\int_{\widetilde{M}}u^{2}},
\]
where the infimum is taken over all smooth functions with compact support. It
is easy to see that for a Cartan-Hadamard manifold $\widetilde{M}^{n}$ with
$K\leq-a^{2}$, we have $\lambda_{0}\left(  \widetilde{M}\right)  \geq\left(
n-1\right)  ^{2}a^{2}/4$.

\section{The Busemann boundary}

Instead of harmonic functions, one can use distance functions to compactify
$\widetilde{M}$. This leads to the Busemann compactification, first introduced
by Gromov in \cite{BGS}. Fix a point $o\in\widetilde{M}$ and define, for
$x\in\widetilde{M}$ the function $\xi_{x}(z)$ on $\widetilde{M} $ by:
\[
\xi_{x}(z)\;=\;d(x,z)-d(x,o).
\]
The assignment $x\mapsto\xi_{x}$ is continuous, one-to-one and takes values in
a relatively compact set of functions for the topology of uniform convergence
on compact subsets of $\widetilde{M}$. The Busemann compactification
$\widehat{M}$ of $\widetilde{M}$ is the closure of $\widetilde{M}$ for that
topology. The space $\widehat{M}$ is a compact separable space. The
\textit{Busemann boundary } $\partial\widehat{M}:=\widehat{M}\setminus
\widetilde{M}$ is made of $1$-Lipschitz continuous functions $\xi$ on
$\widetilde{M}$ such that $\xi(o)=0$ and there exists a sequence $\left\{
a_{k}\right\}  \subset\widetilde{M}$ s.t. $d\left(  o.a_{k}\right)
\rightarrow\infty$ and%
\[
\xi\left(  x\right)  =\lim_{k\rightarrow\infty}d(a_{k},x)-d(a_{k},o),
\]
where the convergence is uniform over compact sets. Elements of $\partial
\widehat{M}$ are called \textit{horofunctions}. We note that this
compactification works for any proper metric space $X$ (cf. \cite{KL}) But in
general, $X$ may fail to be open in its Busemann compactification. This
pathology does not happen for Riemannian manifolds, i.e. we have

\begin{proposition}
$\widetilde{M}$ is open in its Busemann compactification $\widehat{M}$. Hence
the Busemann boundary $\partial\widehat{M}$ is compact.
\end{proposition}

\bigskip For proof see \cite{LW1}. For a Cartan-Hadamard manifold, the
Busemann compactification coincides with the geometric compactification. More precisely,

\begin{proposition}
Let $\widetilde{M}$ be a Cartan-Hadamard manifold and $\left\{  a_{k}\right\}
$ a sequence in $\widetilde{M}$ s.t. $d\left(  o.a_{k}\right)  \rightarrow
\infty$. Let $\sigma_{k}$ be the unique geodesic ray from $o$ to $a_{k}$. Then
$\xi_{a_{k}}$ converges to a horofunction $\xi$ iff $\sigma_{k}$ converges to
a ray $\sigma$. Furthermore, we have
\[
\xi\left(  x\right)  =\lim_{t\rightarrow\infty}d\left(  \sigma\left(
t\right)  ,x\right)  -t.
\]

\end{proposition}

For proof see Ballmann \cite{B} (p30).

Recently, the Busemann compactification has found to be very useful in various
questions, cf. \cite{KL, L1, LW1}. We first describe the application in
\cite{LW1}. Suppose $M$ is a compact Riemannian manifold and $\widetilde{M}$
its universal covering (noncompact). Let $G$ be the fundamental group of $M$
acting on $\widetilde{M}$ isometrically. Observe that we may extend by
continuity the action of $G$ from $\widetilde{M}$ to $\widehat{M}$, in such a
way that for $\xi$ in $\widehat{M}$ and $g$ in $G$,
\[
g.\xi(z)\;=\;\xi(g^{-1}z)-\xi(g^{-1}o).
\]
The volume entropy of $M$ is defined to be the limit
\[
v\left(  g\right)  =\lim_{r\rightarrow\infty}\frac{\ln\mathrm{vol}%
B_{\widetilde{M}}\left(  x,r\right)  }{r},
\]
where $B_{\widetilde{M}}\left(  x,r\right)  $ is the ball of radius $r$
centered at $x$ in the universal covering space $\widetilde{M}$. This
important invariant was introduced by Manning \cite{Ma} who proved

\begin{enumerate}
\item the limit exists and is independent of the center $x\in$ $\widetilde{M}$,

\item $v\leq H$, the topological entropy of the geodesic flow on the unit
tangent bundle of $M$,

\item $v=H$ if $M$ is nonpositively curved.
\end{enumerate}

In \cite{LW1}, we extended the classical theory of Patterson-Sullivan measure
by constructing a family of measures on the \textit{Busemann boundary }
$\partial\widehat{M}$.

\begin{theorem}
There exists a family $\left\{  \nu_{x}:x\in\widetilde{M}\right\}  $of finite
measures on the \textit{Busemann boundary } $\partial\widehat{M}$ s.t.

\begin{enumerate}
\item for any pair $x,y\in\widetilde{M}$ the two measures $\nu_{x}$ and
$\nu_{y}$ are equivalent with%
\[
\frac{d\nu_{x}}{d\nu_{y}}\left(  \xi\right)  =e^{-v\left(  \xi\left(
x\right)  -\xi\left(  y\right)  \right)  },
\]

\item for any $g\in G$ and $x\in\widetilde{M}$%
\[
g_{\ast}\nu_{x}=\nu_{gx}.
\]

\end{enumerate}
\end{theorem}

This family of measures plays a crucial role in the proofs of the following
rigidity results involving the volume entropy.

\begin{theorem}
Let $M^{n}$ be a compact Riemannian manifold with $Ric\geq-\left(  n-1\right)
$. Then the volume entropy satisfies $v\leq n-1$ and equality holds iff $M$ is hyperbolic.
\end{theorem}

\begin{remark}
This result was proved by Knieper \cite{Kn} under the additional assumption
that $M$ is negatively curved.
\end{remark}

As a corollary, in view of the well-known inequality $\lambda_{0}\left(
\widetilde{M}\right)  \leq v^{2}/4$, we deduce the following result which was
previously proved in \cite{W} by a different method.

\begin{theorem}
Let $M^{n}$ be a compact Riemannian manifold with $Ric\geq-\left(  n-1\right)
$. Then $\lambda_{0}\left(  \widetilde{M}\right)  \leq\left(  n-1\right)
^{2}/4$ and equality holds iff $M$ is hyperbolic.
\end{theorem}

\begin{theorem}
Let $M$ be a compact K\"{a}hler manifold with $\dim_{\mathbb{C}}M=m$. If the
bisectional curvature $K_{\mathbb{C}}\geq-2$, then the volume entropy
satisfies $v\leq2m$. Moreover equality holds iff $M$ is complex hyperbolic
(normalized to have constant holomorphic sectional curvature $-4$).
\end{theorem}

\begin{theorem}
Let $M$ be a compact quaternionic K\"{a}hler manifold of $\dim=4m$ with
$m\geq2$ and scalar curvature $-16m\left(  m+2\right)  $. Then the volume
entropy satisfies $v\leq2\left(  2m+1\right)  $. Moreover equality holds iff
$M$ is quaternionic hyperbolic.
\end{theorem}

We refer to the original paper \cite{LW1} for details. More recently, we can
prove some pinching theorems using our method. The first step is the following
rigidity result for $C^{1,\alpha}$ metrics.

\begin{theorem}
Let $M^{n}$ be a (smooth) compact smooth manifold and $g$ a $C^{1,\alpha}$
Riemannian metric. Suppose that $g_{i}$ is a sequence of smooth Riemannian
metrics on $M$ s.t.

\begin{enumerate}
\item $\mathrm{Ric}\left(  g_{i}\right)  \geq-\left(  n-1\right)  $ for each
$i$,

\item $g_{i}\rightarrow g$ in $C^{1,\alpha}$ norm as $i\rightarrow\infty$,

\item the volume entropy $v\left(  g_{i}\right)  \rightarrow n-1$ as
$i\rightarrow\infty$.
\end{enumerate}

Then $g$ is hyperbolic.
\end{theorem}

From this result, we then deduce the following pinching theorem.

\begin{theorem}
There exists a positive constant $\varepsilon=\varepsilon\left(  n,D\right)  $
s.t. if $\left(  M^{n},g\right)  $ is a compact Riemannian manifold of
dimension $n$ satisfying the following conditions

\begin{itemize}
\item $g$ has \textbf{negative} sectional curvature,

\item $\mathrm{Ric}\left(  g\right)  \geq-\left(  n-1\right)  $,

\item $\mathrm{diam}\left(  M,g\right)  \leq D$,

\item the volume entropy $v\left(  g\right)  \geq n-1-\varepsilon$,
\end{itemize}

then $M$ is diffeomorphic to a hyperbolic manifold $\left(  X,g_{0}\right)  $.
\ Moreover, the Gromov-Hausdorff distance $d_{GH}\left(  M,X\right)  $
$\leq\alpha\left(  \varepsilon\right)  \rightarrow0$ as $\varepsilon
\rightarrow0$.
\end{theorem}

We note that this theorem was established by Courtois \cite{C} (unpublished)
in 2000 using the Cheeger-Colding theory. Our proof is different and simpler.
The details will appear in \cite{LW2}.

In \cite{L1} Ledrappier studied another fundamental invariant: the linear
drift (introduced by Guivarc'h \cite{Gu}) which is the following limit for
almost every path $\omega$ of the Brownian motion on $\widetilde{M}$%
\[
l=\lim_{t\rightarrow\infty}\frac{1}{t}d\left(  \omega\left(  0\right)
,\omega\left(  t\right)  \right)  .
\]
If $M$ is negatively curved, Kaimanovich \cite{K1} established a remarkable
integral formula for $l$. Let $\partial\widetilde{M}$ be the geometric
boundary of $\widetilde{M}$. As usual we fix a base point $o\in\widetilde{M}
$. Recall that there is a homeomorphism $\Phi$ from $\partial\widetilde{M}$ to
the Martin boundary $\partial_{\Delta}\widetilde{M}$ by the theorem of
Anderson-Schoen. For each $\xi\in\partial\widetilde{M}$, $h_{\xi}=\Phi\left(
\xi\right)  $ is the unique harmonic function on $\widetilde{M}$ s.t. $h_{\xi
}\left(  o\right)  =1$ and $h_{\xi}\in C\left(  \widetilde{M}^{\ast}%
\backslash\left\{  \xi\right\}  \right)  $ with boundary value zero. With
these notations, the Kaimanovich formula can be written as
\[
l=-\int_{M}\left(  \int_{\partial\widetilde{M}}\left\langle \nabla B_{\xi
},\nabla\ln h_{\xi}\right\rangle \left(  x\right)  h_{\xi}\left(  x\right)
d\nu\left(  \xi\right)  \right)  dm\left(  x\right)  ,
\]
where $\nu$ is the harmonic measure on $\partial\widetilde{M}$ defined by
(\ref{defnu}) and $m$ is the normalized Lebesgue measure on $M$. The main
result in \cite{L1} is a similar integral formula for $l$ in the general case.
The key step is to construct certain measures on the \textit{Busemann boundary
}$\partial\widehat{M}$.

\section{A comparison theorem}

In this section we discuss why compactifications are useful. The basic
principle is that often times a geometric object is much simpler near
infinity. When we look at it on the boundary, \ we capture its essential
features while all the background noise dies off. The first illustration of
this principle is perhaps Mostow's rigidity theorem \cite{Mo}: If
$f:M\rightarrow N$ is a (smooth) homotopy equivalence between two compact
hyperbolic manifolds then $f$ is homotopic to an isometry. In the proof Mostow
considers the lifting $\widetilde{f}:\widetilde{M}\rightarrow\widetilde{N}$
between the universal coverings (which are both $\mathbb{H}^{n}$ in this case)
which is a quasi-isometry. Then $\widetilde{f}$extends to a homeomorphism
$\overline{f}:\partial\widetilde{M}\rightarrow\partial\widetilde{N}$ between
the boundaries. Using the theory of quasi-conformal maps and the fact that the
fundamental group acts on $\partial\widetilde{M}$ ergodically, Mostow shows
that $\overline{f}$ is in fact a Mobius transformation.

More recently, in a seminal paper \cite{BCG1}, Besson, Courtois and Gallot
proved the following theorem which implies the Mostow rigidity theorem in the
rank one cases.

\begin{theorem}
\label{BCG}Let $\left(  N^{n},g_{0}\right)  $ be a compact locally symmetric
space of negative curvature. Let $M^{n}$ be another compact manifold and
$f:M\rightarrow N$ is a continuous map of nonzero degree. Then for any metric
$g$ on $M$

\begin{enumerate}
\item $v\left(  g\right)  ^{n}\mathrm{vol}\left(  M,g\right)  \geq\left\vert
\deg f\right\vert v\left(  g_{0}\right)  ^{n}\mathrm{vol}\left(
N,g_{0}\right)  $;

\item the equality holds iff $f$ is homotopic to an covering map.
\end{enumerate}
\end{theorem}

The proof involves embedding $\widetilde{M}$ into a Hilbert space and a
calibration argument. In a later paper \cite{BCG2}, the same authors gave a
very elegant and simpler proof of their theorem under the additional
assumption that $g$ is also negatively curved. In this second approach, the
Patterson-Sullivan measure on the geometric boundary $\partial\widetilde{M}$
plays a fundamental role. Their method is geometric and flexible and we will
apply it in a slight different situation.

Let $\left(  M^{n},g\right)  $ be a compact Riemannian manifold and
$\pi:\widetilde{M}\rightarrow M$ its universal covering. We pick a base point
$o\in\widetilde{M}$. Let $\partial\widetilde{M}$ be the Martin boundary.

\begin{definition}
For any $p>0$ let%
\[
\beta_{p}\left(  g\right)  =\int_{M}\left(  \int_{\partial\widetilde{M}%
}\left\vert \nabla\log\xi\left(  x\right)  \right\vert ^{2}\xi\left(
x\right)  d\nu\left(  \xi\right)  \right)  ^{p}dm\left(  x\right)  .
\]
Similarly we can consider%
\[
\widetilde{\beta}_{p}\left(  g\right)  =\int_{M}\int_{\partial\widetilde{M}%
}\left\vert \nabla\log\xi\left(  x\right)  \right\vert ^{2p}\xi\left(
x\right)  d\nu\left(  \xi\right)  dm\left(  x\right)  .
\]

\end{definition}

We have $\beta_{p}\left(  g\right)  \leq\widetilde{\beta}_{p}\left(  g\right)
$ by the H\"{o}lder inequality. When $p=1$, $\beta_{p}\left(  g\right)
=\widetilde{\beta}_{p}\left(  g\right)  $ is the Kaimanovich entropy, an
invariant of fundamental importance. This was introduced by Kaimanovich
\cite{K1}. We summarize its main properties:

\begin{enumerate}
\item \cite{K1} $\beta_{1}=\lim_{t\rightarrow\infty}-\frac{1}{t}%
\int_{\widetilde{M}}p_{t}\left(  x,y\right)  \log p_{t}\left(  x,y\right)  dy$
for any $x\in\widetilde{M}$;

\item \cite{K1} $\beta_{1}>0$ iff $\widetilde{M}$ has nonconstant bounded
harmonic functions.

\item \cite{L1, L2} $4\lambda_{0}\leq\beta_{1}\leq v^{2}$, where $v$ is the
volume entropy and $\lambda_{0}$ is the bottom of the $L^{2}$ spectrum of
$\widetilde{M}$.
\end{enumerate}

\bigskip

Let $\left(  N^{n},g_{0}\right)  $ be a compact locally symmetric space of
negative curvature. Theorem \ref{BCG} says that among all metrics $g$ on $N$
with $\mathrm{vol}\left(  N,g\right)  =\mathrm{vol}\left(  N,g_{0}\right)  $
the metric $g_{0}$ has the smallest volume entropy. A natural question is
whether the same is true for the Kaimonovich entropy.

\begin{problem}
\bigskip Let $\left(  N^{n},g_{0}\right)  $ be a compact locally symmetric
space of negative curvature. Is it true that for any metric $g$
\[
\beta_{1}\left(  g\right)  ^{n/2}\mathrm{vol}\left(  N,g\right)  \geq\beta
_{1}\left(  g_{0}\right)  ^{n/2}\mathrm{vol}\left(  N,g_{0}\right)  ?
\]

\end{problem}

We do not know the answer to this question. What we can prove is the following
result which gives an affirmative answer to the same question for $\beta
_{n/2}$.

\begin{theorem}
\label{n/2}Let $\left(  N^{n},g_{0}\right)  $ be a compact locally symmetric
space of negative curvature. Let $M^{n}$ be another compact manifold and
$f:M\rightarrow N$ a (smooth) homotopy equivalence. Then for any metric $g$ on
$M$

\begin{enumerate}
\item $\beta_{n/2}\left(  g\right)  \mathrm{vol}\left(  M,g\right)  \geq
\beta_{n/2}\left(  g_{0}\right)  \mathrm{vol}\left(  N,g_{0}\right)  $;

\item the equality holds iff $f$ is homotopic to an isometry.
\end{enumerate}
\end{theorem}

\bigskip

As a consequence we have the following result which is also an easy corollary
of Theorem \ref{BCG}.

\begin{corollary}
If $g_{0}$ is real hyperbolic and $g$ satisfies $Ric\geq-\left(  n-1\right)
$, then $\mathrm{vol}\left(  M,g\right)  \geq\mathrm{vol}\left(
N,g_{0}\right)  $. Moreover, equality holds iff $g$ is also hyperbolic.
\end{corollary}

This follows from the sharp gradient estimate.

\begin{proposition}
(Li-J. Wang \cite[Lemma 2.1]{LiW}) Let $N^{n}$ be a complete manifold with
$Ric\geq-\left(  n-1\right)  $. If $u$ is a positive harmonic function on $N$,
then $|\nabla\log u|\leq n-1$.
\end{proposition}

We now prove Theorem \ref{n/2}. The homotopy equivalence $f:\left(
M,g\right)  \rightarrow\left(  N,g_{0}\right)  $ induces an isomorphism
$\rho:\Gamma:=\pi_{1}\left(  M\right)  \rightarrow\pi_{1}\left(  N\right)  $.
We view the fundamental groups as groups of deck transformations acting on the
universal covering manifolds. Lifting $f$ we obtain a smooth map
$\widetilde{f}:\widetilde{M}\rightarrow\widetilde{N}$ which is $\Gamma
$-equivariant$\,$, i.e. for any $\gamma\in\Gamma$%
\[
\widetilde{f}\left(  \gamma\cdot x\right)  =\rho\left(  \gamma\right)
\cdot\widetilde{f}\left(  x\right)  .
\]
This is a quasi-isometry and hence $\widetilde{M}$ is Gromov hyperbolic as
$\widetilde{N}$ is. Hence $\widetilde{f}$ extends to a homeomorphism
$\overline{f}:\partial\widetilde{M}\rightarrow\partial\widetilde{N}$ between
the boundaries. By a theorem of Brooks \cite{Br}, $\lambda_{0}\left(
\widetilde{M}\right)  >0$ since $\lambda_{0}\left(  \widetilde{N}\right)  >0$
and the two fundamental groups are isomorphic. Therefore, by Theorem
\ref{Ancona} $\partial\widetilde{M}$ is also the Martin boundary of
$\widetilde{M}$ and let $\left\{  \nu^{x}:x\in\widetilde{M}\right\}  $ be the
harmonic measures. We now define a new map $\widetilde{F}:\widetilde
{M}\rightarrow\widetilde{N}$ applying the construction in \cite{BCG2}:
$\widetilde{F}\left(  x\right)  $ is the barcenter of the measure
$\overline{f}_{\ast}\nu^{x}$ on $\partial\widetilde{N}$, i.e. $\widetilde
{F}\left(  x\right)  $ is the unique minimum point of the following function
on $\widetilde{N}$%
\[
y\rightarrow\int_{\partial\widetilde{N}}B_{\theta}\left(  y\right)  d\left(
\overline{f}_{\ast}\nu^{x}\right)  \left(  \theta\right)  ,
\]
where $B_{\theta}$ is the Busemann function on $\widetilde{N}$ associated to
$\theta\in\partial\widetilde{N}$. For detailed discussion of the barcenter see
\cite{BCG2}. Note that this map is well defined as the support $\nu^{x}$
always has more that two points. By the implicit function theorem, it is easy
to show that $\widetilde{F}$ is smooth. Moreover it is $\Gamma$-equivariant
and hence yields a smooth map $F:M\rightarrow N$. What remains is to estimate
the Jacobian of this map.

By the definition of $F$ we have%

\begin{align*}
& \int_{\partial\widetilde{N}}dB_{\theta}\left(  \widetilde{F}\left(
x\right)  \right)  \left(  \cdot\right)  d\left(  \overline{f}_{\ast}\nu
_{x}\right)  \left(  \theta\right) \\
& =\int_{\partial\widetilde{M}}dB_{\overline{f}\left(  \xi\right)  }\left(
\widetilde{F}\left(  x\right)  \right)  \left(  \cdot\right)  \xi\left(
x\right)  d\nu\left(  \xi\right) \\
& =0.
\end{align*}
Differentiating in $x$ we get
\begin{align*}
& \int_{\partial\widetilde{M}}D^{2}B_{\overline{f}\left(  \xi\right)  }\left(
\widetilde{F}\left(  x\right)  \right)  \left(  \widetilde{F}_{\ast}\left(
x\right)  \left(  \cdot\right)  ,\cdot\right)  \xi\left(  x\right)
d\nu\left(  \xi\right) \\
& =-\int_{\partial\widetilde{M}}dB_{\overline{f}\left(  \xi\right)  }\left(
\widetilde{F}\left(  x\right)  \right)  \left(  \cdot\right)  d\xi\left(
\cdot\right)  d\nu\left(  \xi\right)  .
\end{align*}
We define the following quadratic form $K$ and $H$ on $T_{F\left(  x\right)
}\widetilde{N}$%
\begin{align*}
g_{0}\left(  K_{\widetilde{F}\left(  x\right)  }\left(  u\right)  ,u\right)
& =\int_{\partial\widetilde{M}}D^{2}B_{\theta}\left(  \widetilde{F}\left(
x\right)  \right)  \left(  u,u\right)  \xi\left(  x\right)  d\nu\left(
\xi\right)  ,\\
g_{0}\left(  H_{\widetilde{F}\left(  x\right)  }\left(  u\right)  ,u\right)
& =\int_{\partial\widetilde{M}}dB_{f\left(  \xi\right)  }\left(  \widetilde
{F}\left(  x\right)  \right)  \left(  u\right)  ^{2}\xi\left(  x\right)
d\nu\left(  \xi\right)  .
\end{align*}
Then for any $v\in T_{x}\widetilde{M},u\in T_{\widetilde{F}\left(  x\right)
}\widetilde{N}$
\begin{align*}
& \left\vert g_{0}\left(  K_{\widetilde{F}\left(  x\right)  }\left(  F_{\ast
}\left(  x\right)  \left(  v\right)  \right)  ,u\right)  \right\vert \\
& \leq g_{0}\left(  H_{\widetilde{F}\left(  x\right)  }\left(  u\right)
,u\right)  ^{1/2}\left(  \int_{\partial\widetilde{M}}\frac{\left\vert
\left\langle \nabla\xi\left(  x\right)  ,v\right\rangle \right\vert ^{2}}%
{\xi\left(  x\right)  }d\nu\left(  \xi\right)  \right)  ^{1/2}\\
& =g_{0}\left(  H_{\widetilde{F}\left(  x\right)  }\left(  u\right)
,u\right)  ^{1/2}\left(  \int_{\partial\widetilde{M}}\left\vert \left\langle
\nabla\log\xi\left(  x\right)  ,v\right\rangle \right\vert ^{2}\xi\left(
x\right)  d\nu\left(  \xi\right)  \right)  ^{1/2}.
\end{align*}
Therefore%
\[
\left\vert \det K\right\vert \mathrm{Jac}\widetilde{F}\left(  x\right)
\leq\frac{1}{n^{n/2}}\left\vert \det H\right\vert ^{1/2}\left(  \int
_{\partial\widetilde{M}}\left\vert \nabla\log\xi\left(  x\right)  \right\vert
^{2}\xi\left(  x\right)  d\nu\left(  \xi\right)  \right)  ^{n/2}.
\]
By \cite[Appendix B]{BCG1} we obtain%
\[
\mathrm{Jac}\widetilde{F}\left(  x\right)  \leq\frac{1}{\left(  n+d-2\right)
^{n}}\left(  \int_{\partial\widetilde{M}}\left\vert \nabla\log\xi\left(
x\right)  \right\vert ^{2}\xi\left(  x\right)  d\nu\left(  \xi\right)
\right)  ^{n/2},
\]
where $d=1,2,4$ or $8$ when $\left(  \widetilde{N},g_{0}\right)  $ is the
real, complex, quaternionic hyperbolic space or the Cayley hyperbolic plane,
respectively. Integrating over $M$ yields%
\[
\mathrm{vol}\left(  N,g_{0}\right)  \leq\frac{\mathrm{vol}\left(  M,g\right)
}{\left(  n+d-2\right)  ^{n}}\int_{M}\left(  \int_{\partial\widetilde{M}%
}\left\vert \nabla\log\xi\left(  x\right)  \right\vert ^{2}\xi\left(
x\right)  d\nu\left(  \xi\right)  \right)  ^{n/2}dm\left(  x\right)  .
\]
This proves the inequality as $\beta_{n/2}\left(  g_{0}\right)  =\left(
n+d-2\right)  ^{n}$. If equality holds, it is easy to see that $F$ has to be
an isometry up to a scaling.

\end{document}